\documentclass[12pt]{article}
\def\finproof{\hfill\hbox{\vrule width1.0ex height1.5ex}\vspace{2mm}}

\begin{document}

\begin{center}
{\Large On a class of PDEs with nonlinear distributed in space and time
state-dependent delay term}

\bigskip

{\sc Alexander V. Rezounenko} 
\bigskip

Department of Mechanics and Mathematics, Kharkov University,

 4, Svobody Sqr., Kharkov, 61077, Ukraine

\end{center}
\bigskip

{\bf Abstract.} A new class of nonlinear partial differential
equations with distributed in space and time state-dependent delay
is investigated. We find appropriate assumptions on the kernel
function which represents the state-dependent delay and discuss
advantages of this class. Local and long-time asymptotic
properties, including the existence of global attractor, are
studied.

\bigskip



\bigskip
{\bf 1. Introduction} 
\medskip

Theory of delay differential equations is one of the oldest and
simultaneously, intensively developing branches of the theory of
infinite-dimensional dynamical systems. This theory covers ordinary and
partial delay differential equations, includes studies of discrete and
distributed, finite and infinite delays. Classical methods of differential
equations, theory of distributions and functional analysis allow one to
study wide classes of ordinary and partial differential equations with
delay. We mention only several monographs which are classical references
for delay equations \cite{Hale_book, Walther_book,Azbelev,Mishkis,Wu_book}
and also works which are close to this investigation \cite{travis_webb,
Chueshov-JSM-1992,Cras-1995,Rezounenko-MAG-1997,NA-1998,
Rezounenko-Wu-2006,Rezounenko-JMAA-2007}. Nevertheless, each
(nonlinear) equation requires a separate and careful studying.

Recently, a new class of delay equations attracts attention of
many researchers. These equations have delay (delay term) which
may change, according to the state of the system i.e.
state-dependent (state-selective) delay. The study of such
equations was started in the case of ordinary equations
\cite{Nussbaum-Mallet-1992,Nussbaum-Mallet-1996,Walther_JDE-2003}
and it was recently continued for P.D.E.s in
\cite{Rezounenko-Wu-2006,Rezounenko-JMAA-2007}. For more detailed
discussion and references on delay equations see e.g. introduction
in \cite{Rezounenko-Wu-2006}. We continue our previous
research\cite{Rezounenko-Wu-2006,Rezounenko-JMAA-2007} and present
a wider class of nonlinear equations with distributed in space and
time state-dependent delay terms.

Let us illustrate the main question studied in this article on the
simplified object which is a local in space delay term. Consider
the following simple distributed in time delay term $\int^0_{-r}
b(u(t+\theta, x)) \xi (\theta) d\theta.$ Here function $\xi$
belongs to some space of real valued functions defined on the
delay interval $(-r,0).$ This (kernel) function $\xi$ represents
the rule how the information on the previous stages of the system
(function $u$) is used to model the process. As discussed (see
e.g. \cite{Rezounenko-Wu-2006,Rezounenko-JMAA-2007}), this rule
may change according to the state of the system. Let us denote by
$v\in H$ the state coordinate, where $H$ represents the phase
space. With this notations the state-dependent (state-selective)
delay rule reads $\xi(\theta,v) : (-r,0)\times H \to R$ and the
corresponding delay term becomes $\int^0_{-r} b(u(t+\theta, x))
\xi (\theta,v) d\theta.$ This is a simplified (local) example of
delay terms studied in
\cite{Rezounenko-Wu-2006,Rezounenko-JMAA-2007}. As we will see
(section~3), studying some questions (e.g. stationary solutions),
there is a need to use a wider class of functions $\xi$ (delay
rules) which are space-dependent i.e. $\xi (\theta,x,v).$ For
example, considering a biological system, where $u(t,x)$
represents the density of a population at time moment $t$ at point
$x\in \Omega,$ the delay rule $\xi (\theta,v)$ is the same for all
points $x$ in the domain $\Omega,$ while the delay rule $\xi
(\theta,x,v)$ depends on the point $x\in\Omega$ (e.g. due to the
dependence of food resources on points in $\Omega$). This
interpretation shows that the delay rule $\xi (\theta,x,v)$ is
more realistic biologically and as we will see (section~3)
mathematically. Taking into account the above motivation we need
to find an appropriate class of functions $\xi,$ concentrating on
the character of dependence of $\xi$ on the coordinate
$x\in\Omega.$ This is the main goal of the article. It is
interesting to mention that in spite of the fact that for any time
moment $t\ge 0$ solutions $u(t)$ belong to the space $L^2(\Omega)$
and the phase coordinate $(u(t); u(t+\theta))$ belongs to
$L^2(\Omega)\times L^2(-r,0; L^2(\Omega))$, the values of function
$\xi$ (as functions of $x$) do not necessary belong to
$L^2(\Omega)$. They belong to a wider space $D(A^{-{1/2}})\supset
L^2(\Omega)$ (for more details see theorems~1 and 2 below).

The proposed model has an essential advantage in comparison with the
previous ones (see \cite{Rezounenko-Wu-2006, Rezounenko-JMAA-2007}) to
cover the case of finite and even infinite sequences of isolated
stationary solutions. We also present an algorithm to construct such
state-dependent delay terms.

The article is organized as follows. In section~2 we present the model,
prove the existence and uniqueness of weak solutions, construct the
dynamical system and prove the existence of a global attractor. Section~3
is devoted to stationary solutions and the possibility to use our system
to construct a dynamical system with an a-priory given set of isolated
stationary solutions.
 The results may be applied to the diffusive  Nicholson's blowflies
equation.

\bigskip {\bf 2. Formulation of the model with distributed delay}
\medskip

Consider the following non-local partial differential equation with {\it
state-dependent distributed in space and time delay}
\begin{equation}\label{sdd5-g}
\begin{array}{lll}
&\quad \frac{\partial }{\partial t}u(t,x)+Au(t,x)+du(t,x)\\
&= 
\int^0_{-r} \left\{ \int_\Omega b(u(t+\theta, y))
f(x-y) dy \right\} \xi(\theta,x, u(t), u_t) d\theta \\
&\equiv \big( F(u_t) \big)(x), x\in \Omega ,\end{array}
\end{equation}
 where $A$ is a densely-defined
self-adjoint positive linear operator
 with domain $D(A)\subset L^2(\Omega )$ and with compact
  resolvent, so $A: D(A)\to L^2(\Omega )$ generates an analytic semigroup,
  $\Omega $ is a smooth bounded domain in $R^{n_0}$, $f: \Omega -\Omega \to R$
  is a bounded function to
  be specified later, $b:R\to R$ is a locally Lipschitz bounded map
  ($|b(w)|\le C_b$ with $C_b\ge 0),$ $d$ is a positive constant.
  As usually for delay systems
(see \cite{Hale_book}) for any function $u(t), t\in [a,b], b>a+r$ with
values in a Banach space $X$, we denote by $u_t\equiv u_t(\theta)\equiv
u(t+\theta),$ which is a function of $\theta\in [-r,0]$ with parameter
$t\in [a+r,b].$ Constant $r>0$ is the (maximal) delay of the system.

The function $\xi (\cdot,\cdot,\cdot): [-r,0]\times\Omega\times
H
 \to R$ represents the
state-dependent {\it distributed} delay.
We denote for short $H\equiv L^2(\Omega)\times L^2(-r,0;L^2(\Omega))$ and
also use $\Vert\cdot\Vert$ and $\langle\cdot,\cdot \rangle$ to denote the
norm and scalar product in $L^2(\Omega).$

 We consider equations
(\ref{sdd5-g}) with the following initial conditions
\begin{equation}\label{sdd5-ic}
  u(0+)=u^0\in L^2(\Omega), \quad  u|_{(-r,0)}=\varphi \in L^2 (-r,0;L^2(\Omega)).
\end{equation}
So we write $(u^0,\varphi)\in H.$

\medskip
Now we study the existence and properties of solutions for
distributed delay problem (\ref{sdd5-g}), (\ref{sdd5-ic}).
\medskip

\noindent {\bf Definition 1.} {\it A function $u$ is a {\it weak solution}
of problem (\ref{sdd5-g}) subject to the initial conditions
(\ref{sdd5-ic}) on an interval $[0,T]$ if $u\in L^\infty
(0,T;L^2(\Omega))\cap
 L^2 (-r,T;L^2(\Omega))\cap  L^2 (0,T;D(A^{1\over 2}))$, $u(\theta)=\varphi (\theta)$ for $\theta\in (-r,0)$
 and
\begin{equation}\label{sdd5-sol}
  -\int^T_0\langle u, \dot v\rangle dt +
 \int^T_0\langle A^{1\over 2}u,A^{1\over 2}v\rangle dt +
\int^T_0\langle du-F(u_t),v\rangle dt = -\langle u^0,v(0)\rangle
\end{equation}
 for any function $v\in L^2 (0,T;D(A^{1\over 2}))$ with
 $\dot v\in L^2 (0,T;D(A^{-{1\over 2}}))$ and $v(T)=0.$
  }

\medskip

\bigskip \noindent {\bf Theorem 1.} {\it Assume that
\begin{itemize}
\item[(i)] $b: R\to R$ is locally Lipschitz and bounded i.e.,
there
exists a constant $C_b$ so that that $|b(w)|\le C_b$ 
for all $w\in R$; \item[(ii)] $f: \overline{\Omega-\Omega} \to R$ is
bounded $( |f(\cdot)|\le M_f )$; \item[(iii)] $\xi : [-r,0]\times
\Omega\times L^2(\Omega)\times L^2(-r,0;L^2(\Omega))\to R$ satisfies the
following conditions:\\
{\bf a)} for any $M>0$ there exists $L_{\xi , M}$ so that for all $(v^i,
\psi^i)\in H$ satisfying
 $||v^i||^2+ \int^0_{-r} ||\psi^i(s)||^2 ds\le M^2, i=1,2$
 one has
$$
\quad \int^0_{-r}||\xi (\theta,\cdot, v^1, \psi^1)-\xi (\theta,\cdot, v^2, \psi^2)
 ||_{D(A^{-{1/ 2}})}\,  d\theta$$
\begin{equation}\label{sdd5-xi}\le L_{\xi,M}
  \cdot  \left[ ||v^1-v^2||^2+ \int^0_{-r} ||\psi^1(s)-\psi^2(s)||^2ds
  \right]^{1/2},
\end{equation}
{\bf b)} there exists $C_{(\xi,-{1/ 2})} >0$ so that
\begin{equation}\label{sdd5-cxi}
\int^0_{-r}\Vert \xi(\theta,\cdot\, , v, \psi)\Vert_{D(A^{-{1/ 2}})}\,
d\theta \le C_{(\xi,-{1/ 2})} \,\,\mbox{for all}\,\,  (v, \psi)\in H.
\end{equation}
\end{itemize}

Then for any $(u^0, \varphi) \in H\equiv L^2(\Omega)\times L^2
(-r,0;L^2(\Omega))$ the problem (\ref{sdd5-g}) subject to the initial
conditions (\ref{sdd5-ic}) has a weak solution $u(t)$ on every given time
interval $[0,T]$
 and this solution satisfies
\begin{equation}\label{sdd5-contin} u(t)\in C([0,T];L^2(\Omega)).
\end{equation}
}

\medskip
{\bf Remark.} {\it Properties (iii)-a) and (iii)-b) mean that $\xi$ as a
function of the third and fourth coordinate $(v,\psi)\in H$ is a
(nonlinear) locally Lipschitz and globally bounded mapping $\xi : H \to
L^1 (-r,0;D(A^{-{1\over 2}})).$}

\medskip
\noindent {\it Proof of Theorem~1}. Let us denote by $\{
e_k\}^\infty_{k=1}$ an orthonormal basis of $L^2(\Omega)$ such that
$Ae_k=\lambda_ke_k$, $0< \lambda_1<\ldots<\lambda_k\to +\infty$. We say
that function $u^m(t,x)=\sum\limits^m_{k=1}g_{k,m}(t)e_k(x) $ is a {\it
Galerkin approximate solution of order $m$ for the problem
(\ref{sdd5-g}),(\ref{sdd5-ic})} if
\begin{equation}\label{sdd5-6}
\left\{ \begin{array}{ll} &\langle \dot u^m+Au^m +du^m-F(u^m_{t}),
e_k\rangle =0,\\
&\langle u^m(0+),e_k\rangle =\langle u^0,e_k\rangle , \,\, \langle
u^m(\theta),e_k\rangle=\langle \varphi (\theta) , e_k\rangle
,\,\,\forall \theta\in (-r,0) \end{array}\right.
 \end{equation}
$\forall k=1,\ldots,m$. Here $g_{k,m}\in C^1(0,T;R)\cap
L^2(-r,T;R)$ with $\dot g_{k,m}(t)$ being absolutely continuous.

Equations (\ref{sdd5-6}) for fixed $m$ can be rewritten as a system for
the $m$-dimensional vector-function
$v(t)=v^m(t)=(g_{1,m}(t),\ldots,g_{m,m}(t))^T.$ 
 We notice that $\Vert
u^m(t,\cdot)\Vert^2_{L^2(\Omega)}=
 \sum\limits^m_{k=1}g^2_{k,m}(t)=|v(t)|^2_{R^m}.$

The standard technique (see e.g. \cite{Hale_book}) gives that for any
initial data $\varphi\in L^2(-r,0; R^m),$ $a\in R^m$ there exist
$\alpha>0$ and a unique solution of (\ref{sdd5-6}) $v\in L^2(-r,\alpha;
R^m)$ such that $v_0=\varphi $ and $v(0)=a$, and $v|_{[0,\alpha]}\in
C([0,\alpha]; R^m)$ (for more details see Theorem~6 and Remark~9 from
\cite{Rezounenko-2004} and also Lemma from \cite{Rezounenko-JMAA-2007}).

It is easy to get from (\ref{sdd5-cxi}) and the boundedness of $b$ and $f$
that
$$|\langle F(u_t),v\rangle_{L^2(\Omega)}| =
\left| \int_\Omega \left\{\int^0_{-r} \left[ \int_\Omega b(u(t+\theta, y))
f(x-y) dy\right]  \xi(\theta,x, u(t), u_t) d\theta \right\} v(x) dx
\right| $$ $$= \left| \int^0_{-r} \left[ \int_\Omega b(u(t+\theta, y))
\left\{ \int_\Omega f(x-y) \xi(\theta,x, u(t), u_t)
 v(x) dx \right\} dy\right] d\theta\right|
$$ $$\le C_b M_f |\Omega| \int^0_{-r}
\Vert \xi(\theta,\cdot, u(t), u_t)\Vert_{D(A^{-{1/ 2}})} d\theta \cdot
||v||_{D(A^{{1/ 2}})}.$$ Using (\ref{sdd5-cxi}), one has
\begin{equation}\label{sdd5-F1}
|\langle F(u_t),v\rangle_{L^2(\Omega)}| \le  C_b M_f|\Omega | C_{(\xi,-{1/
2})}\cdot \Vert A^{1/ 2} v\Vert.
\end{equation}


Now, we will get an {\it a-priori} estimate for the Galerkin approximate
solutions for the problem (\ref{sdd5-g}),(\ref{sdd5-ic}). We multiply
(\ref{sdd5-6}) by $g_{k,m}$ and sum over $k=1,\cdots ,m$. Hence for
$u(t)=u^m(t)$ and $t\in (0,\alpha]\equiv(0,\alpha (m)]$, the local
existence interval for $u^m(t)$, we get
\begin{equation}\label{sdd5-11}
{1\over 2}{d\over dt}\Vert u(t)\Vert^2 + \Vert A^{1/2}u(t)\Vert^2 + d
\Vert u(t)\Vert^2  \le |\langle F(u_t),u(t)\rangle|.
\end{equation}
Using (\ref{sdd5-F1}), (\ref{sdd5-11}) we obtain
\begin{equation}\label{sdd5-12}
{d\over dt}\Vert u(t)\Vert^2 + \Vert A^{1/2}u(t)\Vert^2 + 2d\Vert
u(t)\Vert^2\le C^2_b M^2_f|\Omega |^2 C^2_{(\xi,-{1/ 2})}\equiv\tilde k_1.
\end{equation}
Since ${d\over dt}\Vert u(t)\Vert^2 + \Vert A^{1/2}u(t)\Vert^2 +2d\Vert
u(t)\Vert^2= {d\over dt}\left(\Vert u(t)\Vert^2 + \int^t_0\Vert
A^{1/2}u(\tau)\Vert^2d\tau \right.$ $\left. +2d\int^t_0\Vert
u(\tau)\Vert^2d\tau\right) $, we denote by $\chi (t)\equiv \Vert
u(t)\Vert^2 + \int^t_0\Vert A^{1/2}u(\tau)\Vert^2d\tau +2d\int^t_0\Vert
u(\tau)\Vert^2d\tau$ and rewrite the last estimate as follows $ {d\over
dt}\chi (t) \le
\tilde k_1.$
We obtain 
$ \chi(t)\le \chi(0)+{\tilde k_1 t}=\Vert u(0)\Vert^2 +\tilde k_1 t.$ So,
we have the {\it a -priori} estimate
\begin{equation}\label{sdd5-14}
\Vert u(t)\Vert^2 + \int^t_0\Vert
A^{1/2}u(\tau)\Vert^2d\tau+2d\int^t_0\Vert u(\tau)\Vert^2d\tau\le \Vert
u(0)\Vert^2+\tilde k_1 t.
\end{equation}
 Estimate (\ref{sdd5-14}) gives that, for
$u^0\in L^2(\Omega)$ 
the family of approximate solutions $\{ u^m(t)\}^\infty_{m=1}$ is
uniformly (with respect to $m\in {\bf N}$) bounded in the space $L^\infty
(0,T;L^2(\Omega))\cap L^2 (0,T;D(A^{1/2})),$ where $D(A^{1/2})$ is the
domain of the operator $A^{1/2}$ and $[0,T]$ is the local existence
interval. From (\ref{sdd5-14}) we also get the continuation of $u^m(t)$ on
any interval, so (\ref{sdd5-14}) holds for all $t>0.$

Using the definition of Galerkin approximate solutions (\ref{sdd5-6}) and
their property (\ref{sdd5-14}), we can integrate over $[0,T]$ to obtain
$\int^T_0\Vert A^{-{1/2}} \dot u^m(\tau)\Vert^2 d\tau \le C_T$ for any
$T.$ These properties of the family $\{ u^m(t)\}^\infty_{m=1}$ give that
$\{ (u^m(t); \dot u^m(t))
\}^\infty_{m=1}$ is a bounded sequence in the space 
\begin{equation}\label{sdd5-13}
X_T\equiv L^\infty (0,T;L^2(\Omega))\cap L^2
(0,T;D(A^{1/2}))\times L^2(0,T;D(A^{-{1/2}})).
\end{equation}
 Then
there exist a function $(u(t); \dot u(t))$ and a subsequence $\{
u^{m_k} \} \subset \{ u^m \}$ such that
\begin{equation}\label{sdd5-*}
(u^{m_k}; \dot u^{m_k}) \quad \hbox{*-weakly converges to}
 \quad  (u; \dot u) \quad \hbox{in the space} \quad X_T.
\end{equation}
By a standard argument (using the strong convergence $u^{m_k}\to u$ in the
space $L^2 (0,T;L^2(\Omega))$ which follows from (\ref{sdd5-*}) and the
Doubinskii's theorem, one can show (see e.g. Lions (1969), Chueshov (1999)
and Rezounenko (1997)) that any *-weak limit is a solution of
(\ref{sdd5-g}) subject to the initial conditions (\ref{sdd5-ic}). To prove
the continuity of
weak solutions we use the well-known (see also \cite[thm. 1.3.1]{Lions-Magenes-book}) 

\medskip
\noindent {\bf Proposition~1} (Proposition 1.2 in
\cite{showalter}). {\it Let the Banach space $V$ be dense and
continuously embedded in the Hilbert space $X;$ identify $X=X^*$
so that $V\hookrightarrow X\hookrightarrow V^*.$ Then the Banach
space $W_p(0,T)\equiv \{ u\in L^p(0,T;V) : \dot u\in
L^q(0,T;V^*)\}$ (here $p^{-1}+q^{-1}=1$) is contained in
$C([0,T];X).$ }

In our case $X=L^2(\Omega), V=D(A^{{1/2}}), V^*=D(A^{-{1/2}}), p=q=1/2$
(see (\ref{sdd5-13}),(\ref{sdd5-*})). Hence Proposition~1 gives
(\ref{sdd5-contin}). The proof of Theorem~1 is complete. \rule{5pt}{5pt}

\bigskip

Now we describe a sufficient condition for the uniqueness of weak
solutions.
\medskip

\noindent {\bf Theorem~2.} {\it Assume that functions $b$ and $f$ are as
in Theorem~1 (satisfy properties (i),(ii)), function $\xi$ satisfies
property (iii)-a) and
\begin{equation}\label{sdd5-23}
   \xi (\cdot,\cdot, v,\psi) \in L^\infty(-r,0;D(A^{-{1/ 2}}))\,\,\mbox{for all}\,\,  (v, \psi)\in H.
\end{equation}
\nopagebreak Then solution of (\ref{sdd5-g}), (\ref{sdd5-ic}) given by
Theorem~1 is unique. }

\medskip

\noindent {\it Proof of Theorem~2.} Let $u^1$ and $u^2$ be two solutions
of (\ref{sdd5-g}), (\ref{sdd5-ic}).
Below we denote for short $w(t)=w^{m}(t)=u^{1,m}(t)-u^{2,m}(t)$ - the
difference of corresponding Galerkin approximate solutions. Hence
\begin{equation}\label{sdd5-24}
{d\over dt}\Vert w(t)\Vert^2 + 2 \Vert A^{1/2}w(t)\Vert^2 + 2d \Vert
w(t)\Vert^2 = \langle F(u^1_t)- F(u^{2}_t),w(t)\rangle.
\end{equation}
Let us consider the difference $\langle F(u^1_t)- F(u^{2}_t),w(t)\rangle$
in details  (see (\ref{sdd5-g})).
$$\langle F(u^1_t)-
F(u^{2}_t),w(t)\rangle\equiv \int_\Omega\left[
 \int^0_{-r} \left\{ \int_\Omega b(u^1(t+\theta, y)) f(x-y) dy
\right\} \xi(\theta,x, u^1(t), u^1_t) d\theta  - \right.$$
$$
-\left. \int^0_{-r} \left\{ \int_\Omega b(u^{2}(t+\theta, y)) f(x-y) dy
\right\} \xi(\theta,x, u^{2}(t), u^{2}_t) d\theta \right] \cdot w(t,x) dx
$$
$$=\int_\Omega\left[ \int^0_{-r} \left\{ \int_\Omega b(u^1(t+\theta,
y)) f(x-y) dy \right\} \xi(\theta,x, u^1(t), u^1_t) d\theta -\right.
$$
$$
-\left. \int^0_{-r} \left\{ \int_\Omega b(u^{2}(t+\theta, y)) f(x-y) dy
\right\} \xi(\theta,x, u^1(t), u^1_t) d\theta \right] \cdot w(t,x) dx,
$$
$$
+\int_\Omega\left[ \int^0_{-r} \left\{ \int_\Omega b(u^{2}(t+\theta, y))
f(x-y) dy \right\} \xi(\theta,x , u^1(t), u^1_t) d\theta -\right.
$$
$$
-\left. \int^0_{-r} \left\{ \int_\Omega b(u^{2}(t+\theta, y)) f(x-y) dy
\right\} \xi(\theta,x,u^{2}(t),u^{2}_t) d\theta \right] \cdot w(t,x) dx.
$$
Using the local Lipschitz property of $b,$ (\ref{sdd5-23}) and
(\ref{sdd5-xi}), we deduce 
$$|\langle F(u^1_t)-F(u^2_t), w(t)\rangle | \le L_b M_f \int^0_{-r}
\left\{ \int_\Omega |w(t+\theta,y)|\, dy \cdot \int_\Omega
|\xi(\theta,x,u^{1}(t),u^{1}_t) | \cdot |w(t,x)|\, dx\right\} d\theta$$
$$ + C_b M_f |\Omega | \int^0_{-r}||\xi(\theta,\cdot,u^{1}(t),u^{1}_t)-
\xi(\theta,\cdot,u^{2}(t),u^{2}_t) ||_{D\left(A^{-{1/ 2}}\right)} d\theta
\cdot ||A^{{1/ 2}}w(t)||$$
$$\le L_b M_f \sqrt{|\Omega |} \int^0_{-r}
||w(t+\theta,\cdot)||\cdot ||\xi(\theta,\cdot,u^{1}(t),u^{1}_t)
||_{D\left(A^{-{1/ 2}}\right)} \cdot ||A^{{1/ 2}}w(t)||\, d\theta$$
$$ + C_b M_f |\Omega | L_{\xi,M}\left[ ||w(t)||^2 + \int^0_{-r} ||w(t+s)||^2 ds\right]^{1/ 2}
\cdot ||A^{{1/ 2}}w(t)||$$
$$\le L_b M_f \sqrt{|\Omega |} \quad \hbox{ess }\hskip-3mm \sup_{\theta\in (-r,0)} ||\xi(\theta,\cdot,u^{1}(t),u^{1}_t)
||_{D\left(A^{-{1/ 2}}\right)} \cdot \int^0_{-r}
||w(t+\theta,\cdot)||\cdot ||A^{{1/ 2}}w(t)||\, d\theta$$
$$ + {1\over 2}\, ||A^{{1/ 2}}w(t)||^2 +{1\over 2}\, C^2_b M^2_f |\Omega |^2
L^2_{\xi,M}\left[ ||w(t)||^2 + \int^0_{-r} ||w(t+s)||^2 ds\right]
$$
$$\le {1\over 2}\, ||A^{{1/ 2}}w(t)||^2 + {1\over 2} \, L^2_b M^2_f
|\Omega |\, r  \left[\hbox{ess }\hskip-3mm \sup_{\theta\in (-r,0)}
||\xi(\theta,\cdot,u^{1}(t),u^{1}_t) ||_{D\left(A^{-{1/ 2}}\right)}
\right]^2 \cdot \int^0_{-r} ||w(t+\theta,\cdot)||^2\, d\theta
$$
$$ + {1\over 2}\, ||A^{{1/ 2}}w(t)||^2 +{1\over 2}\, C^2_b M^2_f |\Omega |^2
L^2_{\xi,M}\left[ ||w(t)||^2 + \int^0_{-r} ||w(t+s)||^2 ds\right].
$$

Finally, we get the existence of positive constants $C_1, C_2$ such that
$$|\langle F(u^1_t)-
F(u^{2}_t),w(t)\rangle | \le ||A^{1/2}w(t)||^2 + C_1\int^0_{-r}
||w(t+\theta)||^2 d\theta + C_2||w(t)||^2  $$


The last estimate and (\ref{sdd5-24}) give 
$${d\over dt}\Vert w(t)\Vert^2 + 2\Vert A^{1/2}w(t)\Vert^2 +2d\Vert w(t)\Vert^2
\le ||A^{1/2}w(t)||^2 + C_1\int^0_{-r} ||w(t+\theta)||^2 d\theta +
C_2||w(t)||^2 $$

$$\le ||A^{1/2}w(t)||^2 + C_1\left( \int^0_{-r}||w(\theta)||^2 d\theta
+\int^t_{0}||w(s)||^2 ds\right) + C_2||w(t)||^2.
$$
Hence
$${d\over dt}\Vert w(t)\Vert^2 + \Vert A^{1/2}w(t)\Vert^2 +2d\Vert w(t)\Vert^2
\le C_1\left( \int^0_{-r}||w(\theta)||^2 d\theta +\int^t_{0}||w(s)||^2
ds\right) + C_2||w(t)||^2
$$
and property $\Vert A^{1/2}v\Vert^2\ge \lambda_1\Vert v\Vert^2$ gives

$${d\over dt}\left[ \Vert w(t)\Vert^2 + (\lambda_1+2d)\int^t_{0}\Vert
w(s)\Vert^2 ds\right] \le C_1\left( \int^0_{-r}||w(\theta)||^2 d\theta
+\int^t_{0}||w(s)||^2 ds\right) + C_2||w(t)||^2.
$$
It implies that there exists $C_3>0,$ such that for $Z(t)\equiv \Vert
w(t)\Vert^2 + (\lambda_1+2d)\int^t_{0}\Vert w(s)\Vert^2 ds,$ we have
$${d\over dt} Z(t) \le C_3Z(t)+ C_1\int^0_{-r}||w(\theta)||^2 d\theta.$$
Gronwall lemma implies
\begin{equation}\label{sdd5-25}
Z(t)\le \left( \Vert w(0)\Vert^2 + C_1C^{-1}_3\int^0_{-r}||w(\theta)||^2
d\theta\right) \cdot e^{C_3t}.
\end{equation}
The last estimate allows one to apply the well-known

\smallskip

\noindent {\bf Proposition~2.} \cite[Theorem~9]{yosida}
 {\it Let $X$ be a Banach space. Then any *-weak convergent sequence
 $\{ w_k\}^\infty_{n=1}\in X^{*}$  *-weak converges to an element
 $w_\infty\in X^{*}$ and $\Vert w_\infty\Vert_X \le\liminf_{n\to\infty} \Vert w_n\Vert_X.$
}

Hence, for the difference $u^1(t)-u^2(t)$ of two solutions we have
$$\Vert u^1(t)-u^2(t)\Vert^2 + 2(\lambda_1+d)\int^t_{0}\Vert
u^1(s)-u^2(s)\Vert^2ds $$
\begin{equation}\label{sdd5-26}
\le \left( \Vert u^1(0)-u^2(0)\Vert^2 +
C_1C^{-1}_3\int^0_{-r}||\varphi^1(\theta)-\varphi^2(\theta)||^2
d\theta\right) \cdot e^{C_3t}.
\end{equation}

We notice that by (\ref{sdd5-contin}) the difference $\Vert
u^1(t)-u^2(t)\Vert$ makes sense for all $t\in [0,T],\, \forall T>0.$ The
last estimate gives the uniqueness of solutions and completes the proof of
Theorem~2. \rule{5pt}{5pt}

\medskip

Theorems~1 and 2 allow us to define the evolution semigroup $S_t : H\to
H$, with  $H\equiv L^2(\Omega)\times L^2 (-r,0;L^2(\Omega))$, by the
formula $S_t(u^0;\varphi)\equiv (u(t);u(t+\theta)), \, \theta\in (-r,0),$
where $u(t)$ is the weak solution of (\ref{sdd5-g}),(\ref{sdd5-ic}). The
continuity of the semigroup with respect to time follows from
(\ref{sdd5-contin}), and with respect to initial conditions from
(\ref{sdd5-26}).

For the study of long-time asymptotic properties of the above
evolution semigroup we recall (see e.g. \cite{Babin-Vishik,
Temam_book})

\medskip

\noindent {\bf Definition~2.} {\it A global attractor of the
semigroup $S_t$ is a closed bounded set ${\mathcal U}$ in $H,$
strictly invariant ($S_t{\mathcal U}={\mathcal U}$ for any $t\ge
0$), such that for any bounded set $B\subset H$ we have
$\lim\limits_{t\to +\infty} \sup \{ dist_H (S_ty, {\mathcal U}),
y\in B\} =0.$ }

\medskip


\noindent {\bf Theorem 3.} {\it Assume functions $b$ and $f$
satisfy properties (i), (ii) of Theorem~1. Let function $\xi$
satisfy  properties (iii)-a) of Theorem~1, (\ref{sdd5-23}) and
also there exists $C_{\xi,0}>0$ such that (c.f. (\ref{sdd5-cxi}))
\begin{equation}\label{sdd5-27}
\int^0_{-r}\Vert \xi(\theta,\cdot\, , v, \psi)\Vert\, d\theta \le
C_{(\xi,0)} \,\,\mbox{for all}\,\,  (v, \psi)\in H.
\end{equation}
Then the dynamical system $(S_t;H)$ has a compact global attractor
${\mathcal U}$ which is a bounded set in the space $H_1\equiv
D(A^{\alpha})\times W$, where $W=\{ \varphi : \varphi \in L^\infty
(-r,0;D(A^{\alpha})), \dot\varphi \in L^\infty (-r,0;D(A^{\alpha-1}))\}$,
$\alpha\le {1\over 2}.$}

\medskip

{\it Proof of Theorem~3.} To prove the existence of the global
attractor we use classical theorem saying that it is sufficient
for the dynamical system $(S_t,H)$ to be dissipative and
asymptotically compact (see \cite{Babin-Vishik, Temam_book,
Chueshov_book}).

The property (\ref{sdd5-27}) gives the estimate stronger than
(\ref{sdd5-F1}):
\begin{equation}\label{sdd5-28}
|\langle F(u_t),v\rangle_{L^2(\Omega)}| \le  C_b M_f|\Omega |
C_{(\xi,0)}\cdot \Vert v\Vert
\end{equation}
which is necessary for the property of dissipativeness of
$(S_t;H).$ The rest of the proof, including the property of
asymptotic compactness, is standard (see e.g. \cite{Babin-Vishik,
Chueshov_book,Rezounenko-MAG-1997} and also
\cite{Rezounenko-Wu-2006,Rezounenko-JMAA-2007})~.~\finproof

\bigskip


\centerline{\bf 3. Stationary solutions} \nopagebreak
\medskip

For simplicity of presentation, in this section we consider operator
$A=(-\Delta_D)>0,$ where $\Delta_D$ is the Laplace operator in
$L^2(\Omega)$ with the Dirichlet boundary conditions. In this case (which
is sufficient for the application to the Nicholson's blowfly equation), we
have $D(A)=H^2(\Omega)\cap H^1_0(\Omega),$ $D(A^{1/ 2})=H^1_0(\Omega),
D(A^{-{1/ 2}})=H^{-1}(\Omega).$ For more details on this classical Sobolev
spaces see e.g. \cite{Lions-Magenes-book}.

In this section we concentrate on the stationary solutions. First of all,
by definition~1 (of a weak solution), $u(t,x)\in L^2(0,T;D(A^{{1/
2}}))=L^2(0,T;H^{1}_0(\Omega)),$ so for the stationary solution
$u(t,x)\equiv u^{st} (x),$ one has $u^{st}\in H^{1}_0(\Omega).$

Let us consider an arbitrary 
function $u^{st}\in H^{1}_0(\Omega)\subset L^2(\Omega).$ Our goal is to
find conditions on a function $\xi(\cdot,\cdot,\cdot,\cdot)$ such that the
system (\ref{sdd5-g}) has stationary solution $u(t)\equiv u^{st}\in
H^{1}_0(\Omega)$ for all ${t\in R}.$ Let us denote by
$\overline{u^{st}}\equiv \overline{u^{st}}(\theta)\equiv u^{st}, \theta\in
[-r,0].$

Since for $u^{st}= 0\in H^{1}_0(\Omega)$ we can choose $\xi (\cdot, \cdot,
0,0)\equiv 0,$ we concentrate below on the case $u^{st}\neq 0\in
H^{1}_0(\Omega).$

From (\ref{sdd5-g}) and $\frac{\partial }{\partial t}u(t,x)\equiv 0,$ we
have

\begin{equation}\label{sdd5-29}
\begin{array}{lll}
&\quad Au^{st}(x)+d\cdot u^{st}(x) =
 \int_\Omega b(u^{st}(y))
f(x-y) dy \cdot \int^0_{-r}\xi(\theta,x, u^{st}, \overline{u^{st}})
d\theta, \quad x\in \Omega. \end{array}
\end{equation}
As we will show, it is sufficient to define in a proper way the value of
$\xi$ for the second and third coordinates equal $(u^{st},
\overline{u^{st}})\in H\equiv L^2(\Omega)\times L^2(-r,0;L^2(\Omega))$
only. We propose to look for this value, decomposing it on the time and
space coordinates i.e.
\begin{equation}\label{sdd5-30}
\xi(\theta,x, u^{st}, \overline{u^{st}}) = \chi(\theta)\cdot \hat v(x),
\quad \theta\in [-r,0],\quad  x\in \Omega.
\end{equation}
Now equation (\ref{sdd5-29}) reads
\begin{equation}\label{sdd5-31}
\begin{array}{lll}
&\quad Au^{st}(x)+d\cdot u^{st}(x) =
 \int_\Omega b(u^{st}(y))
f(x-y) dy \cdot \hat v(x)\cdot \int^0_{-r}\chi(\theta) d\theta, \quad x\in
\Omega. \end{array}
\end{equation}

We need the following elementary

\medskip

{\bf Lemma.} {\it Assume $u^{st}\neq 0\in L^2(\Omega).$ Let the function
$f$ be
strictly positive and $f\in C^\infty (\overline{\Omega-\Omega}).$
Let function $b$ be bounded and satisfy $b(w)>0$ for all $w\neq 0.$
%

Then the function
\begin{equation}\label{sdd5-32}
p(x)\equiv \int_\Omega b(u^{st}(y)) f(x-y) dy
\end{equation}
satisfies properties:
 $p\in C(\overline{\Omega}),$   $\inf \{ p(x) : x \in
\overline{\Omega} \} \equiv p_{min}> 0$ and \ $\sup \{
\frac{\partial}{\partial x_i }p(x) : x \in \overline{\Omega}, i=1,...,n_0
\} \equiv p^\prime_{max}<\infty.$

}

\medskip

The continuity of $p$ on $\overline{\Omega}$ follows immediately from the
continuity of $f$ and the Cauchy-Schwartz inequality
$$
|p(x^1)- p(x^2)| = \left| \int_\Omega b(u^{st}(y)) \left[ f(x^1-y) -
f(x^2-y) \right] dy\right|
$$
\begin{equation}\label{sdd5-33}
\le || b(u^{st}(\cdot ))|| \cdot \left[ \int_\Omega \left|f(x^1-y) -
f(x^2-y)\right|^2 dy \right]^{1/2}.
\end{equation}
We also use that for all $y\in \Omega$ one has $\left|(x^1-y) -
(x^2-y)\right| =\left| x^1 - x^2\right|.$

Properties $b(w)>0$ for all $w\neq 0,$ $u^{st}\neq 0\in L^2(\Omega)$ and
strict positivity of $f$ imply that $p(x)>0$ for all $x \in
\overline{\Omega}.$ Hence, the continuity of $p$ and the Weierstrass
theorem give $\inf \{ p(x) : x \in \overline{\Omega} \} \equiv p_{min} >
0.$ The boundedness of partial derivatives of $p$ is due to $f\in C^\infty
(\overline{\Omega-\Omega}).$
\finproof

We also assume
\begin{equation}\label{sdd5-34}
\int^0_{-r}\chi(\theta) d\theta\neq 0.
\end{equation}

Under the assumptions of Lemma and (\ref{sdd5-34}) we have $p(x)\cdot
\int^0_{-r}\chi(\theta) d\theta\neq 0$ for all $x \in \Omega.$ So we can
write (see (\ref{sdd5-31}))
\begin{equation}\label{sdd5-35}
\hat v(x)=\frac{Au^{st}(x)+d\cdot u^{st}(x)}{\int_\Omega b(u^{st}(y))
f(x-y) dy \cdot \int^0_{-r}\chi(\theta) d\theta}.
\end{equation}
We notice that (\ref{sdd5-35}) is the equality in $D(A^{-{1/
2}})=H^{-1}(\Omega)$ (in the sense of distributions). As we saw, by
definition~1 (of a weak solution), $u(t,x)\in L^2(0,T;D(A^{{1/
2}}))=L^2(0,T;H^{1}_0(\Omega)),$ so for the stationary solution
$u(t,x)\equiv u^{st} (x),$ one has $u^{st}\in H^{1}_0(\Omega).$ This
implies $Au^{st}\in H^{-1}(\Omega)=D(A^{-{1/ 2}})$ and $Au^{st}+d\cdot
u^{st} \in H^{-1}(\Omega).$ To show that $\hat v\in H^{-1}(\Omega)$ we
remind the following

\medskip

{\bf Proposition~3.} \cite[Theorem~12.1]{Lions-Magenes-book} {\it Let $m$
be positive integer. Then any element $h\in H^{-m}(\Omega)$ may be
represented (in the non-unique way) in the form
$$ h=\sum_{|j|\le m} D^j h_j, \quad h_j \in L^{2}(\Omega).
$$
}

Here $D^\alpha \equiv \frac{\partial^{\alpha_1+\ldots+\alpha_n}}{\partial
x^{\alpha_1}_1 \ldots \partial x^{\alpha_n}_n},\quad \alpha=\{
\alpha_1,\ldots,\alpha_n\}, \quad |\alpha|=\alpha_1+\ldots +\alpha_n.$

\smallskip

In our case $m=1$ and if we denote by $h=Au^{st}+d\cdot u^{st} \in
H^{-1}(\Omega)$ and by $q(x)\equiv p^{-1}(x),$ then $\hat v\in
H^{-1}(\Omega)$ reads as $qh\in H^{-1}(\Omega)$. Using proposition~3, we
write $h=h_0 + \sum^{n_0}_{i=1} \frac{\partial}{\partial x_i } h_i, \quad
h_i \in L^{2}(\Omega). $ Hence,
\begin{equation}\label{sdd5-36}
\hat v=qh=qh_0 + \sum^{n_0}_{i=1} q\frac{\partial}{\partial x_i } h_i
=qh_0 + \sum^{n_0}_{i=1} \left[ \frac{\partial}{\partial x_i } (qh_i) -
h_i\frac{\partial}{\partial x_i } q \right].
\end{equation}

{\bf Remark. } {\it We notice that all the derivatives are understood in
the sense of distributions (see \cite{schwartz, kolmogorov-fomin}). The
term $q\cdot h$ is understood as the distribution which is obtained by
multiplication of the distribution $h$ by the infinitely differentiable
function $q$ (by definition, $(q\cdot h,\varphi)\equiv ( h,q\cdot\varphi),
\forall \varphi\in D$ as in \cite{schwartz, kolmogorov-fomin}), since the
operation of multiplication is not defined for two distributions. Using
this definition, it is easy to check that $\frac{\partial}{\partial x_i }
(q\cdot h) = h\cdot \frac{\partial}{\partial x_i } q +
q\cdot\frac{\partial}{\partial x_i } h.$ }

By proposition~3, to get $\hat v\in H^{-1}(\Omega)$ it is enough to show
(see (\ref{sdd5-36})) that
\begin{equation}\label{sdd5-37}
q h_i \in L^2(\Omega), \qquad h_i \frac{\partial}{\partial x_i }q  \in
L^2(\Omega). \end{equation}

 The first inclusion in (\ref{sdd5-37}) follows from Lemma and
 $$\int_\Omega |q(x) h_i(x)  |^2 \, dx \le
 \left[ \sup \{ |q(x)| : x\in \overline{\Omega}\}\right]^2
 \cdot ||h_i||^2=p^{-2}_{min}\cdot ||h_i||^2<+\infty.$$
The second inclusion in (\ref{sdd5-37}) holds due to
$$\int_\Omega |h_i(x) \frac{\partial}{\partial x_i }q(x)   |^2 \, dx \le
\left[ \sup \left\{ \left|\frac{\partial}{\partial x_i }q(x)\right| : x\in
\overline{\Omega}, i=1,...,n_0\right\}\right]^2\cdot ||h_i||^2 $$ $$\le
\left[\frac{p^\prime_{max}}{p^2_{min}}\right]^2\cdot ||h_i||^2<+\infty.$$
Here we use (see Lemma)
$$\left|\frac{\partial}{\partial x_i }q(x)\right| = \left|\frac{\partial}{\partial x_i }
p^{-1}(x)\right| = \left|-\frac{\partial p(x)}{\partial x_i }\cdot
p^{-2}(x)\right|
$$
$$\le \sup \left\{ \left|\frac{\partial}{\partial x_i }q(x)\right| : x\in
\overline{\Omega}, i=1,...,n_0\right\} p^{-2}_{min}\le p^\prime_{max}\cdot
p^{-2}_{min}.
$$
So we get the property $\hat v\in H^{-1}(\Omega)=D(A^{-{1/ 2}})$ which is
very important for us to justify the choice of assumptions on the
state-dependent function $\xi$ (the choice of a class of functions $\xi$)
in this article. Now we see that, assuming (in addition to
(\ref{sdd5-34})) that $\int^0_{-r} |\chi (\theta)|\, d\theta <\infty, $
the function  $\xi$ defined by (\ref{sdd5-30}) with $\hat v$ defined by
(\ref{sdd5-35}) possesses the property (\ref{sdd5-cxi}).

As a result, we may conclude that for any (finite or infinite) sequence of
{\it isolated} points  $\{ u^{st,k}\} \subset H^1_0(\Omega)\subset
 L^2(\Omega)$ we can define a state-dependent function $\xi,$ which
satisfies assumptions of Theorem~1 and such that system (\ref{sdd5-g})
with this $\xi$ will have all the points $\{ u^{st,k}\} \subset
H^1_0(\Omega)$ as stationary solutions $ u^k(t)\equiv u^{st,k}, t\in R$.
The last property means that our model with distributed in space and time
state-dependent delay term may be successfully used  having information
(say from experiments) on an arbitrary set of isolated stationary
solutions.

We notice that the definition of values of
$\xi(\cdot,\cdot,u^{st,k},\overline{u^{st,k}})$ by (\ref{sdd5-30}),
(\ref{sdd5-35}) on a set of isolated points does not contradict property
(\ref{sdd5-xi}) since the last one deals with the case of convergent
sequence of points in $H.$

To conclude  this section we collect all the assumptions on functions
$\xi, b, f$ used in our considerations:

\begin{itemize}
\item[{\bf Ab)}] Function $b: R\to R$ is locally Lipschitz, bounded
and satisfies $b(w)>0$ for all $w\neq 0.$
\item[{\bf Af1)}] Function $f: \overline{\Omega-\Omega}\to R$ is bounded.

\item[{\bf Af2)}] Function $f$ is strictly positive and $f\in C^\infty
(\overline{\Omega-\Omega}).$

\item[{\bf A$\xi 1$)}] Function $\xi : [-r,0]\times \Omega\times
L^2(\Omega)\times L^2(-r,0;L^2(\Omega))\to R$ satisfies the
following condition:\\
for any $M>0$ there exists $L_{\xi , M}$ so that for all $(v^i,
\psi^i)\in H$ satisfying
 $||v^i||^2+ \int^0_{-r} ||\psi^i(s)||^2 ds\le M^2, i=1,2$
 one has
$$\quad \int^0_{-r}||\xi (\theta,\cdot, v^1, \psi^1)-\xi (\theta,\cdot, v^2, \psi^2)
||_{D(A^{-{1/ 2}})}\,  d\theta$$ 
$$\le L_{\xi,M}
  \cdot  \left[ ||v^1-v^2||^2+ \int^0_{-r} ||\psi^1(s)-\psi^2(s)||^2ds
  \right]^{1/2}.
$$

\item[{\bf A$\xi 2$)}]
There exists $C_{(\xi,-{1/ 2})} >0$ so that\\
$\int^0_{-r}\Vert \xi(\theta,\cdot\, , v, \psi)\Vert_{D(A^{-{1/ 2}})}\,
d\theta \le C_{(\xi,-{1/ 2})} \,\,\mbox{for all}\,\,  (v, \psi)\in H.
$

\item[{\bf A$\xi 3$)}] Function $\xi$ satisfies\quad
 $  \xi (\cdot,\cdot, v,\psi) \in L^\infty(-r,0;D(A^{-{1/ 2}}))\,\,\mbox{for all}\,\,  (v, \psi)\in H.
$

\item[{\bf A$\xi 4$)}]
There exists $C_{(\xi,0)} >0$ so that
$\int^0_{-r}\Vert \xi(\theta,\cdot\, , v, \psi)\Vert\, d\theta \le
C_{(\xi,0)} \,\,\mbox{for all}\,\,  (v, \psi)\in H.
$

\item[{\bf A$\chi $)}] Function $\chi $ satisfies $\int^0_{-r}\chi(\theta)
d\theta\neq 0$ and $\int^0_{-r} |\chi (\theta)|\, d\theta <\infty.$

\end{itemize}


\medskip

As an application we can consider the diffusive Nicholson's blowflies
equation (see e.g. \cite{So-Yang,So-Wu-Yang}) with state-dependent delays.
More precisely, we consider equation (\ref{sdd5-g}) where $-A$ is the
Laplace operator with the Dirichlet boundary conditions, $\Omega\subset
R^{n_0}$ is a bounded domain with a smooth boundary, the function $f$ can
be a constant as in \cite{So-Yang,So-Wu-Yang} which leads to the local in
space coordinate term or, for example, $ f(s)={1\over \sqrt{4\pi\alpha}}
e^{-s^2/4\alpha}$, as in \cite{So-Wu-Zou} which corresponds to the
non-local term, the nonlinear function $b$ is given by $b(w)=p\cdot
we^{-w}.$ Function $b$ is bounded and $b(w)>0$ for all $w\neq 0.$
As a result, we conclude that for any functions $\xi$ satisfying
conditions of
Theorems~2 and 3 
the dynamical system $(S_t,H)$ has a global attractor (Theorem~3).

So, our system (\ref{sdd5-g}) with distributed in space and time
state-dependent delay term may be successfully used to study Nicholson's
blowflies equation with an arbitrary set of isolated stationary solutions.

\medskip

\noindent {\bf Acknowledgements.}
The author wishes to thank Hans-Otto Walther for bringing state-dependent
delay differential equations to his attention.

\medskip

\hfill November 20, 2006

\end{document}